%% file: beta-normal-cleaned.tex
\RequirePackage[l2tabu, orthodox]{nag}
\RequirePackage{snapshot}

\documentclass[10pt,onecolumn]{extarticle}

\makeatletter
\if@twocolumn
  \usepackage[dvips,letterpaper,top=0.5in, bottom=0.5in, left=0.75in, right=0.5in,includefoot,heightrounded]{geometry}
\else
  \usepackage[dvips,letterpaper,margin=1in,includefoot,heightrounded]{geometry}
\fi

\usepackage{srcltx}

\usepackage[russian,portuges,english]{babel}

\iflanguage{portuges}
    {\newcommand{\keywordname}{Palavras-chaves}}
    {\newcommand{\keywordname}{Keywords}}

\usepackage{amsmath}
\usepackage{amssymb,amsfonts}

\usepackage{abstract}

\usepackage{graphicx}
\usepackage{epsfig}
\usepackage[usenames,dvipsnames,svgnames,x11names]{xcolor}
\usepackage{subfigure}

\usepackage{booktabs}

\usepackage{setspace}
\usepackage{flushend}

\usepackage{cite}

\usepackage{hyperref}
\usepackage[normalem]{ulem}

\usepackage{enumerate}

\usepackage[noend]{algpseudocode}

\usepackage{listings}

\usepackage[short,12hr]{datetime}
       \usepackage{fouriernc}
\newtheorem{proposition}{Proposition}

\def\QED{\mbox{$\square$}}
\def\proof{\noindent{\it Proof:~}}
\def\endproof{\hspace*{\fill}~\QED\par\endtrivlist\unskip}

\makeatletter

\newcommand{\printtitle}{%
\makeatletter
\if@twocolumn

\twocolumn[%
  \maketitle
  \begin{onecolabstract}
    \myabstract
  \end{onecolabstract}
  \begin{center}
    \small
    \textbf{\keywordname}
    \\\medskip
    \mykeywords
  \end{center}
  \bigskip
]
\saythanks
\else
  \maketitle
  \begin{onecolabstract}
    \myabstract
  \end{onecolabstract}
  \begin{center}
    \small
    \textbf{\keywordname}
    \\\medskip
    \mykeywords
  \end{center}
  \bigskip
  \onehalfspacing
\fi
\makeatother
}

\author{%
L. C. R\^ego%
\thanks{
L. C. R\^ego was with the
Departamento de Estat\'{\i}stica,
Universidade Federal de Pernambuco (UFPE), Brazil;
currently he
is with the  Departamento de Estat\'{\i}stica e Matem\'atica Aplicada,
Universidade Federal do Cear\'a, Brazil.
E-mail: \url{leandro@dema.ufc.br}}
\and
R. J. Cintra%
\thanks{%
R. J. Cintra is with the
Signal Processing Group,
UFPE, Brazil.
E-mail: \url{rjdsc@de.ufpe.br}}
\and
G. M. Cordeiro%
\thanks{%
G. M. Cordeiro is with the
Departamento de Estat\'{\i}stica,
UFPE, Brazil.}
}

\title{%
On Some Properties of the Beta Normal Distribution}

\newcommand{\myabstract}{%
The beta normal distribution is a generalization of both the normal
distribution and the normal order statistics. Some of its mathematical
properties and a few applications have been studied in the literature.
We provide a better foundation for some properties and an analytical
study of its bimodality.
The hazard rate function and the limiting behavior
are examined.
We derive explicit expressions for moments,
generating function,
mean deviations using a power series expansion for
the quantile function,
and Shannon entropy.
}

\newcommand{\mykeywords}{%
Beta normal distribution,
bimodality,
hazard function,
generating function,
quantile function,
mean deviation,
Shannon entropy.
}

\date{}

\begin{document}

\printtitle

\section{Introduction}

The beta normal (BN) distribution~\cite{eugene2002beta} contains as special
sub-models the normal distribution and the normal order statistics. In fact, for the $\mathrm{BN}(\alpha,\beta,\mu,\sigma)$
distribution with parameters $\alpha>0$, $\beta>0$,
$\mu\in\mathbb{R}$ and $\sigma>0$, the probability density
function (\mbox{pdf}) is
\begin{align}
\label{fdpbn}
f(x)=
\frac{1}{\sigma\mathrm{B}(\alpha,\beta)}\,
\left[\Phi\left( \frac{x-\mu}{\sigma} \right)
\right]^{\alpha-1}\,
\left[1-\Phi\left(\frac{x-\mu}{\sigma} \right)
\right]^{\beta-1}\,
\phi\left(\frac{x-\mu}{\sigma} \right),
\quad
x\in\mathbb{R},
\end{align}
where $\Phi(\cdot)$ and $\phi(\cdot)$ are the cumulative distribution
function (\mbox{cdf}) and the standard normal density function,
respectively, and $\mathrm{B}(\cdot,\cdot)$ is the beta function.
The parameters $\alpha$ and $\beta$ are shape parameters,
that characterize the skewness, kurtosis and bimodality of (\ref{fdpbn}),
$\mu$ is a location parameter and $\sigma$ is a dispersion parameter
that stretches out or shrinks the distribution. The BN distribution can be both
unimodal and bimodal. A few authors have provided important properties
of the BN distribution~\cite{gupta2004moments}.

Location and scale parameters are redundant in (\ref{fdpbn}) since if
$Z\sim \mathrm{BN}(\alpha,\beta,0,1)$ then $X=\mu+\sigma\,Z\,\sim\mathrm{BN}(\alpha,\beta,\mu,\sigma)$.
The random variable $Z$ has the so-called beta standard normal (BSN) distribution.
The parameters $\alpha$ and $\beta$ control skewness through the relative tail weights.
The BN distribution
is symmetric if $\alpha=\beta$;
presents negative skewness when $\alpha<\beta$,
and positive skewness when $\alpha>\beta$.
Moreover,
as $\beta$ decreases, the skewness increases.
Conversely,
as $\alpha$ increases, the skewness increases.

For $\alpha>1$ and $\beta>1$, the BN
distribution has positive excess kurtosis, and as $\alpha$ and $\beta$ increase
the higher its peak.
On the other hand,
when $\alpha<1$ and $\beta< 1$, it has negative excess kurtosis,
and as both $\alpha$ and $\beta$ decrease until bimodality, the heavier are the tails.

As pointed out in~\cite{famoye2004bimodality},
the bimodal distribution can occur
in fields as diverse as:
statistical processing of
synthetic aperture radar imagery~\cite{el-zaart2007sar},
neurological disorder assessment~\cite{louisa2007age},
distribution of radiation mechanisms~\cite{nardini2008afterglows},
and quantification of atmospheric pressure~\cite{Zangvil2001pressure}.

Indeed, a careful choice of parameters allows the BN distribution to exhibit bimodality.
Famoye and collaborators pioneered the study of the shape of this distribution~\cite{famoye2004bimodality}.
Although their study addressed several issues about the bimodality
properties, it was primarily numerical in nature, relying on estimation techniques, for instance.
In particular, the study on the parameter range for which the BN distribution switches
behavior from unimodality to bimodality was partial.
It was asserted that ``\emph{the BN distribution becomes bimodal for certain values
of the parameters $\alpha$ and $\beta$ and the analytical solution
of $\alpha$ and $\beta$, where the distribution becomes bimodal, cannot
be solved algebraically}''~\cite{famoye2004bimodality}.

We examine some structural properties of the BN distribution that were not
developed so far. The rest of the paper is organized as follows.
In Section~2, we provide an analytical study of the bimodality region of the BN distribution.
We propose an exact algebraic description for the critical bimodal parameter values.
In Section~3, we investigate the hazard rate function and its limiting behavior.
Explicit expressions for the moments are discussed in Section~4.
Further, in Sections~5-8, we derive quantile measures, generating function,
mean deviations and Shannon entropy.

\section{Bimodality}

The analysis of the critical points of the BN density function furnishes
a natural path for characterizing the distribution shape and
quantifying the number of modes.
After considering the normalization
$z=\frac{x-\mu}{\sigma}$,
we have
\begin{align*}
\frac{\partial f(z)}{\partial z}=
&
\frac{1}{\mathrm{B}(\alpha,\beta)}
\phi(z)
[\Phi(z)]^{\alpha-2}
[1-\Phi(z)]^{\beta-2}\times\\
&
\Big\{
(\alpha-1)\phi(z)[1-\Phi(z)]
-
(\beta-1)\phi(z)\Phi(z)
-
z \Phi(z) [1-\Phi(z)]
\Big\}.
\end{align*}
Let us refer to the term in curly brackets as $s(z)$.
At the critical points, where $\frac{\partial f(z)}{\partial z}=0$,
we have $s(z)=0$, since the remaining terms of $\frac{\partial f(z)}{\partial z}$ are
strictly positive. Hence, the critical points satisfy the following implicit equation
\begin{align}
z=(2-\alpha-\beta)\frac{\phi(z)}{[1-\Phi(z)]}+
(\alpha-1)\frac{\phi(z)}{\Phi(z)[1-\Phi(z)]}.
\label{criticalpoints}
\end{align}

\begin{figure*}
\centering
\subfigure[]{\epsfig{file=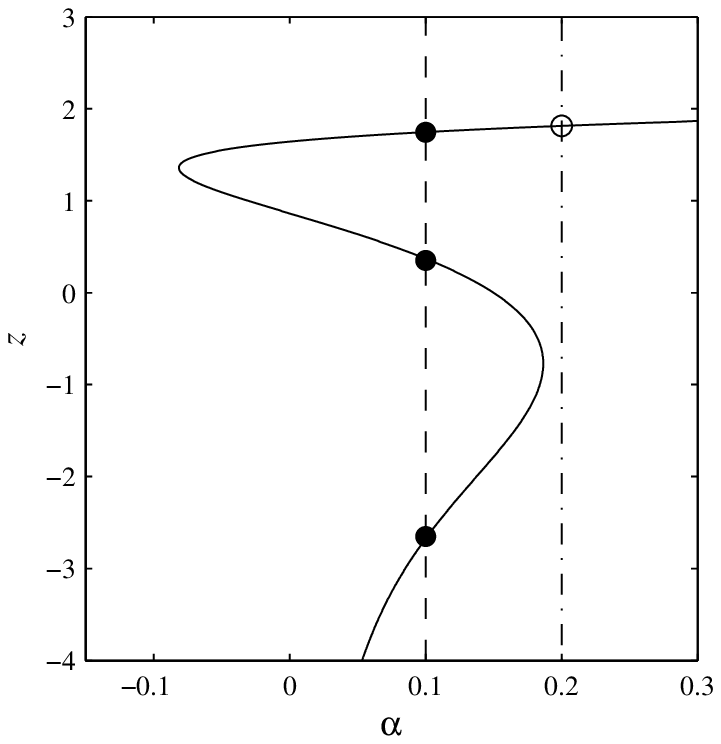}}
\qquad
\subfigure[]{\epsfig{file=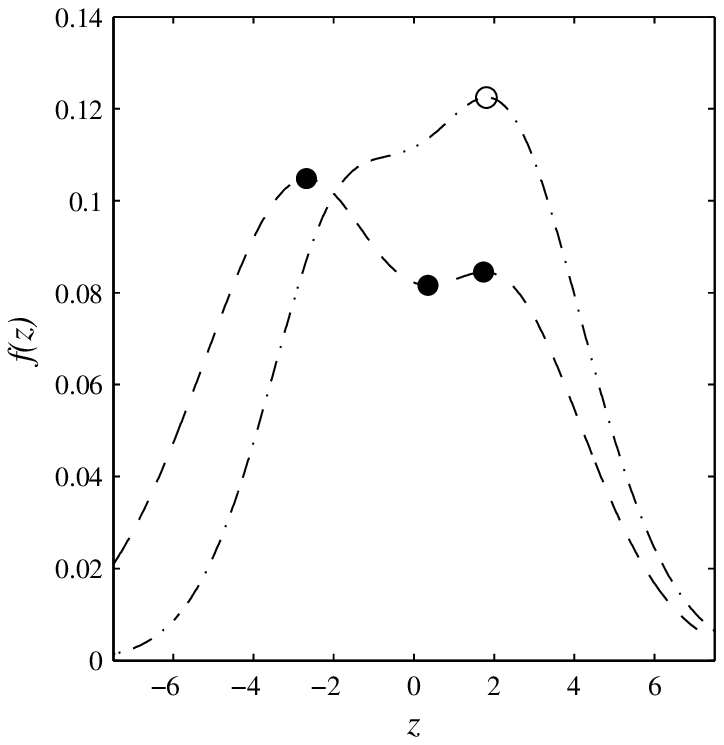}}
\caption{
(a)~Implict curve for the critical values of $f(z)$ for a fixed value of $\beta=0.15$.
(b)~Associated density plots for $\alpha=0.1$ (dashed) and $\alpha=0.2$ (dot-dashed).
}
\label{fig1}
\end{figure*}

In~\cite{famoye2004bimodality}, it is claimed that the solutions
of (\ref{criticalpoints}) are the modes of the distribution.
In fact, there are solutions of this equation that can be
local minimum of the density function, thus they are not modes. For example,
Figure~\ref{fig1}(a) gives the plot of the solutions of (\ref{criticalpoints})
in terms of $\alpha$ for a fixed value of $\beta=0.15$.
For $\alpha=0.1$, there are three solutions indicated by filled dots ($\bullet$).
However, for $\alpha=0.2$, only one critical point is found ($\circ$).

Figure~\ref{fig1}(b) provides the associated density plots.
For $\alpha=0.1$, only two of the marked points are indeed modes of the distribution
but the remaining point characterizes a local minimum. For $\alpha=0.2$, the single
mode is depicted. For given $\beta$, the choice of $\alpha$ determines the number of modes of $f(z)$.
Additionally, we note that the only parts of the implicit curve
with probabilistic meaning are those situated in the region $\alpha>0$. Conversely, the above
discussion also holds if $\alpha$ is fixed and $\beta$ varies.

In order to determine which critical points are modes of the distribution,
we should consider the sign of the second derivative at the critical points.
In particular, a mode of $f(z)$ is a critical point with non-positive second derivative.
At the critical points, we have
\begin{align*}
\frac{\partial^2 f(z)}{\partial z^2}=&
\frac{1}{\mathrm{B}(\alpha,\beta)}\,
\phi(z)\,[\Phi(z)]^{\alpha-2}\,[1-\Phi(z)]^{\beta-2}\,
\frac{\partial s(z)}{\partial z}.
\end{align*}
The sign of $\frac{\partial^2 f(z)}{\partial z^2}$ is the
same of $\frac{\partial s(z)}{\partial z}$.
Then, at a mode, the condition $\frac{\partial s(z)}{\partial z}\leq 0$ holds.
Explicit evaluation of $\frac{\partial s(z)}{\partial z}$ yields
\begin{align}
\frac{\partial s(z)}{\partial z}=-\Phi(z)\,[1-\Phi(z)]- z\,\phi(z)\,
\left[\alpha-(\alpha+\beta)\,\Phi(z)\right]+(2-\alpha-\beta)\,\phi^2(z).
\label{derivativeofs(z)}
\end{align}

We consider the variational behavior of the critical points of $f(z)$ with respect
to changes in the parameter $\alpha$. From (\ref{criticalpoints}), the first derivative
of $z$ with respect to $\alpha$ is
\begin{align*}
\frac{\partial z}{\partial \alpha}
=
&
\frac{\phi(z)\,[1-\Phi(z)]}{\Phi(z)\,[1-\Phi(z)]+ z\,\phi(z)
\left[\alpha-(\alpha + \beta)\,\Phi(z)\right]
-(2-\alpha-\beta)\,\phi^2(z)}.
\end{align*}

Since $\phi(z)[1-\Phi(z)]>0$, the sign of $\frac{\partial z}{\partial \alpha}$
depends entirely on the behavior of the denominator term.
Moreover, except for the sign,
this denominator is equal to $\frac{\partial s(z)}{\partial z}$.
Thus, for any $z$, we have
\begin{align*}
\frac{\partial z}{\partial \alpha}&=\frac{\phi(z)[1-\Phi(z)]}%
{-\frac{\partial s(z)}{\partial z}}
\quad
\text{and}
\quad
\mathrm{sign}\left(\frac{\partial z}{\partial \alpha}\right)
=
-\mathrm{sign}\left(\frac{\partial s(z)}{\partial z}\right),
\end{align*}
where $\mathrm{sign}(\cdot)$ is the sign function.
Further, if $\frac{\partial s(z)}{\partial z}$ is negative at a critical point $z$,
then $z$ must be a mode which is an increasing function of $\alpha$. However, nothing
prevents the existence of a mode for which $\frac{\partial s(z)}{\partial z}$ vanishes.
Indeed, in the case $\alpha=\beta =1-\pi/4$, it can be shown that
\begin{align*}
\frac{\partial f(0)}{\partial z}
=
\frac{\partial^2 f(0)}{\partial z^2}
=
\frac{\partial^3 f(0)}{\partial z^3}
=
0
\quad
\text{and}
\quad
\frac{\partial^4 f(0)}{\partial z^4}
&=
\frac{4\sqrt{2}}{\pi}
\frac{\Gamma(3/2-\pi/4)}{\Gamma(1-\pi/4)}
\left(\frac{3}{\pi}-1 \right)<0.
\end{align*}

\begin{figure}
\centering
\epsfig{file=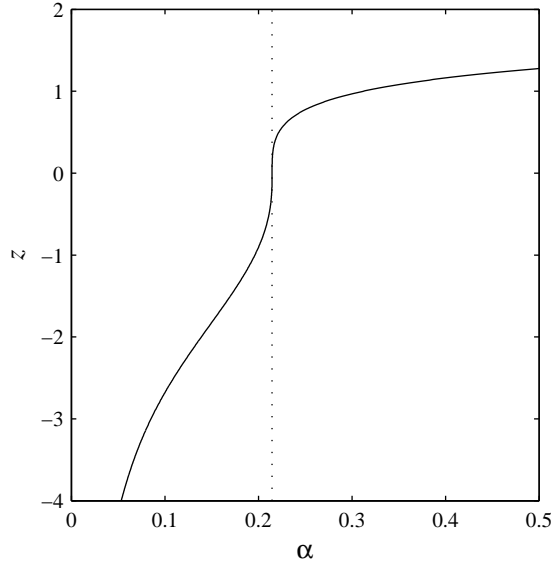}
\caption{At $\alpha=\beta=1-\pi/4$, the quantity $\partial z/\partial \alpha$ is undefined.}
\label{fig.inflexao}
\end{figure}

Thus, $z=0$ is a mode of the distribution.
Numerical computations give evidence that
this is the unique case for which the second derivative vanishes
at the mode.
In this situation,
$\frac{\partial z}{\partial \alpha}$ is undefined.
Nevertheless, it is still true that $z$ is an
increasing function of $\alpha$, as shown in Figure~\ref{fig.inflexao}.
We are then in position to state the following proposition.

\begin{proposition}
If $z$ is a mode location,
then $z$ is an increasing function of $\alpha$.
\end{proposition}
Based on a similar analysis, we can show that every mode is a decreasing function of $\beta$.
These results were previously examined in~\cite[Corollary 3]{famoye2004bimodality}.
However, their proof relied on an inaccurate derivation of
$\frac{\partial z}{\partial \alpha}$
and
$\frac{\partial z}{\partial \beta}$,
which were mistakenly shown to be strictly positive and negative,
respectively.

Now, we consider the symmetric case when $\alpha=\beta$. The following proposition was also
stated in~\cite[Corollary 1]{famoye2004bimodality}.
However, the associated proof was partial.
Indeed, it was only shown that if $z_0$ is a critical point,
then $-z_0$ is also a critical point of the distribution;
we complete the proof here.

\begin{proposition}
Let $\alpha=\beta$. If $z_0$ is a modal point, then $-z_0$ is also a modal point.\label{propsymmetric}
\end{proposition}

\proof
For $\alpha=\beta$, $f(z)$ is an even function and, consequently,
$\frac{\mathrm{d}f(z)}{\mathrm{d}z}$ is odd.
Thus, if $z_0$ is a critical point of $f(z)$, so is $-z_0$.
Note that $s(z)$ is also odd which implies
that $\frac{\partial s(z)}{\partial z}$ is even.
Then, $\frac{\partial s(z_0)}{\partial z}=\frac{\partial s(-z_0)}{\partial z}$,
which assures that $-z_0$ is a modal point.
\endproof

We also provide a more complete proof of the following result stated
in~\cite[Corollary 2]{famoye2004bimodality}.

\begin{proposition}
If $\mathrm{BN}(\alpha,\beta,\mu,\sigma)$ has a mode at $z_0$,
then $\mathrm{BN}(\beta,\alpha,\mu,\sigma)$ has a mode at $-z_0$.
\end{proposition}

\proof
If $\alpha=\beta$, then the result follows from Proposition~\ref{propsymmetric}.
We address the case $\alpha\ne\beta$. Let $t(z)$ be the result of interchanging the roles
of $\alpha$ and $\beta$ in $s(z)$. Note that $t(-z)=-s(z)$ and $z_0$ is
a critical point of $\mathrm{BN}(\beta,\alpha,\mu,\sigma)$
if, and only if, $t(z_0)=0$. Thus,
if $z_0$ is a critical point of $\mathrm{BN}(\alpha,\beta,\mu,\sigma)$,
then $-z_0$ is a critical point of $\mathrm{BN}(\beta,\alpha,\mu,\sigma)$.
Moreover,
since $\alpha\ne\beta$,
it follows that $z_0$ is a mode of $\mathrm{BN}(\beta,\alpha,\mu,\sigma)$
if, and only if,
$\frac{\partial t(z_0)}{\partial z}<0$.
Since
$\frac{\partial t(z_0)}{\partial z}=\frac{\partial s(-z_0)}{\partial z}$,
the result follows.
\endproof

\begin{figure}
\centering
\epsfig{file=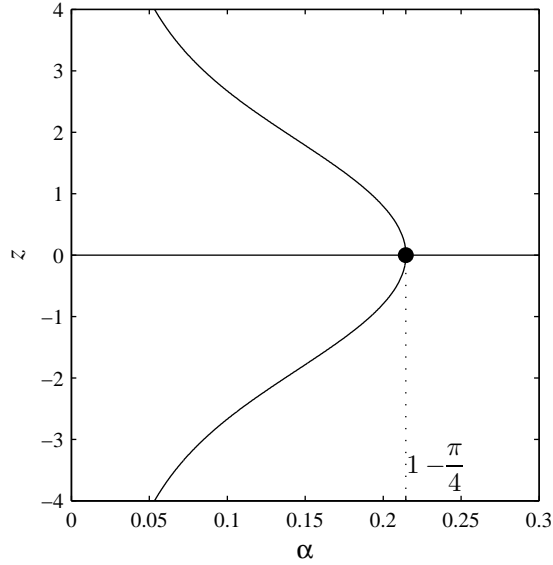}
\caption{Implict curve for the critical points for $\alpha=\beta$.
The locus $z=0$ is part of the curve.}
\label{fig.alpha=beta}
\end{figure}

Consider again the symmetric case $\alpha=\beta$. In Figure~\ref{fig.alpha=beta},
we plot the critical points of $f(z)$ in terms of the parameter $\alpha$.
This curve is obtained from (\ref{criticalpoints}).
There is a critical value of $\alpha$
after which the BN density function exhibits a single critical point ($z=0$),
that is the unique mode of the distribution.
This particular critical value can be determined as follows.
From~(\ref{criticalpoints}), we can express $\alpha$ in terms of $z$ as
\begin{align*}
\alpha=1+z\frac{\Phi(z)}{\phi(z)} \frac{1-\Phi(z)}{1-2\Phi(z)},
\quad \text{for $z\ne 0$}.
\end{align*}

From Figure~\ref{fig.alpha=beta}, we note that the exact critical value $\alpha^*$
can be obtained as the limit of the above function as $z$ tends to zero.
Thus, using L'H\^opital rule, we obtain
\begin{align*}
\alpha^*&=
1+\lim_{z\rightarrow 0}z\frac{\Phi(z)}{\phi(z)}\frac{1-\Phi(z)}{1-2\Phi(z)}
=
1+\frac{\sqrt{2\pi}}{4}\lim_{z\rightarrow 0}\frac{z}{1-2\Phi(z)}
\\
&
=
1+\frac{\sqrt{2\pi}}{4}\lim_{z\rightarrow 0}\frac{1}{-2\phi(z)}
=
1-\frac{\pi}{4}
\approx
0.2146.
\end{align*}

Further, in the symmetric case, we have $s(0)=0$.
Thus, $z=0$ is always a critical point.
We examine the sign of $\frac{\partial^2 f(z)}{\partial z^2}$ at $z=0$.
Since its sign is the same as that of $\frac{\partial s(z)}{\partial z}$,
by~(\ref{derivativeofs(z)}), we obtain
\begin{align*}
\left.\frac{\partial s(z)}{\partial z}\right|_{z=0}
=
-\frac{1}{4}+\frac{1-\alpha}{\pi}.
\end{align*}
Hence,
$\left.\frac{\partial s(z)}{\partial z}\right|_{z=0}<0$
if, and only if,
$\alpha>1-\frac{\pi}{4}$.
As discussed before,
in the symmetric case,
where $\alpha=\beta=1-\pi/4$,
$z=0$ is the unique mode of the distribution.

\subsection{Modality Regions}

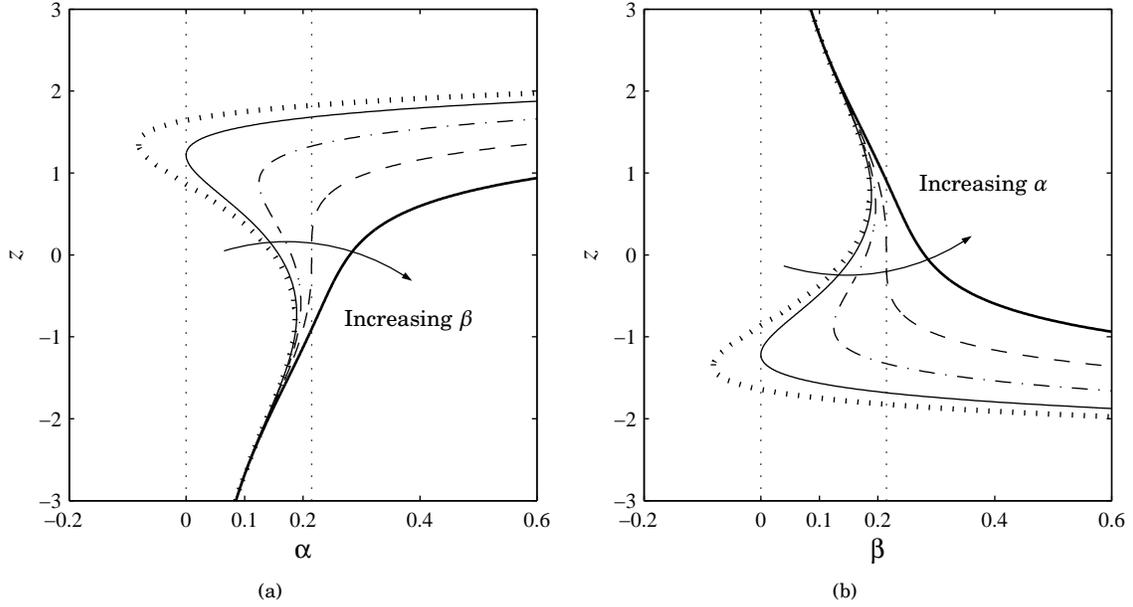
\begin{figure*}
\centering
\subfigure[]{\input{figure4a.pstex_t}}
\subfigure[]{\input{figure4b.pstex_t}}
\caption{
Implict curve for the critical points of $f(z)$
(a) as a function of $\alpha$
for $\beta \in \{0.15, 0.158896, 0.18, 1-\pi/4, 0.28\}$
and
(b) as a function of $\beta$
for $\alpha \in \{0.15, 0.158896, 0.18, 1-\pi/4, 0.28\}$.
}
\label{fig.several.curves}
\end{figure*}

In Figure~\ref{fig.several.curves}(a), we plot the critical points $z$ of the
distribution in terms of $\alpha$ for selected values of $\beta$.
Concerning the modality of the BN distribution,
we can separate three regions:
(i)~$0<\beta\leq 0.158896$,
(ii)~$0.158896< \beta < 1 - \pi/4$,
and
(iii)~$\beta \geq 1-\pi/4$.

For $0<\beta\leq 0.158896$ and small positive values of $\alpha$,
the distribution is already bimodal,
but there is a critical value of $\alpha$ after which
the distribution has a single critical point and then it is unimodal.
For $0.158896<\beta<1-\pi/4$ and small positive values of $\alpha$,
the distribution is unimodal and then
there is a critical value of $\alpha$ after which
the distribution has three critical points,
where two of them are modes of the distribution.
But, as in the previous case,
when the value of $\alpha$ increases,
there is another critical value of $\alpha$ after which the
distribution has a single critical point and then it is unimodal.
For the last case,
$\beta \geq 1-\pi/4$ and the distribution is always unimodal.

From (\ref{criticalpoints}),
we can express $\alpha$ in terms of $z$ for a fixed value
of $\beta$.
Let $\alpha_\beta(z)$ denote this function of $z$ by fixing $\beta$.
Hence,
\begin{align*}
\alpha_\beta(z)=
z\frac{\Phi(z)}{\phi(z)} - \frac{(2-\beta)\Phi(z)-1}{1-\Phi(z)}.
\end{align*}

Moreover, for a fixed value of $\beta\in (0,1-\pi/4)$,
let $\alpha^*_\beta$ denote the local maximum of $\alpha_\beta(z)$.

A similar analysis of the variation of the critical points
in terms of $\beta$ for several values of $\alpha$ can be made.
Figure~\ref{fig.several.curves}(b) illustrates this analysis.
Due to its symmetric behavior,
the discussion follows \emph{mutatis mutandis}.
So,
from (\ref{criticalpoints}),
we can express $\beta$ in terms of $z$ for a fixed
value of $\alpha$.
Let $\beta_\alpha(z)$ denote such function given by
\begin{align*}
\beta_\alpha(z)=
\frac{\alpha-1}{\Phi(z)} - z \frac{1-\Phi(z)}{\phi(z)} + 2 - \alpha.
\end{align*}
Moreover,
for a fixed value of $\alpha\in (0,1-\pi/4)$,
let $\beta^*_\alpha$ denote the local maximum of $\beta_\alpha(z)$.
Using the symmetry properties of $\phi(\cdot)$ and $\Phi(\cdot)$,
we obtain
\begin{align*}
\beta_\alpha(-z)
=
\frac{\alpha-1}{1-\Phi(z)} + z \frac{\Phi(z)}{\phi(z)} + 2 - \alpha.
\end{align*}
After simple manipulations, for a fixed quantity $\gamma$,
we have $\beta_\gamma(-z) = \alpha_\gamma(z)$. For any real function $f(z)$,
the set of values of the local maxima of $f(z)$ and $f(-z)$
are exactly the same. Therefore,
for a given $0<\gamma<1-\pi/4$, we have $\beta^*_\gamma=\alpha^*_\gamma$.
Table~\ref{table1} lists the numerical values of
$\beta^*_\gamma=\alpha^*_\gamma$
for several values of $\gamma\in(0,1-\pi/4)$.
Thus, the region of the parameters $\alpha$ and $\beta$ for which
the distribution is bimodal reduces to
\begin{align*}
\left\{
(\alpha,\beta):
0<\alpha<\min\left(1-\frac{\pi}{4},\alpha^*_\beta\right),
0<\beta<\min\left(1-\frac{\pi}{4},\beta^*_\alpha\right)
\right\}
.
\end{align*}
The curves
$\{(\alpha,\beta^*_\alpha):0<\alpha<1-\pi/4\}$
and
$\{(\alpha^*_\beta,\beta):0<\beta<1-\pi/4\}$
which delimit the above region are
symmetric with respect to the line $\alpha=\beta$.

\begin{figure}
\centering
\epsfig{file=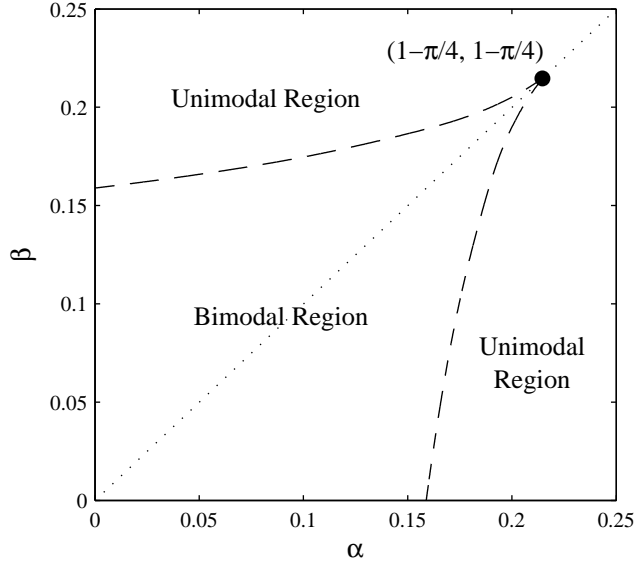}
\caption{Modality regions.}
\label{fig.boundary}
\end{figure}

In~\cite{famoye2004bimodality},
a numerical estimate of these curves was approximated
by means of linear regression techniques.
However, such reported curves were not symmetric
with respect to the line $\alpha=\beta$.
More accurately, Figure~\ref{fig.boundary}
shows the modality regions of the BN distribution.

\begin{table}
\centering
\caption{Selected boundary coordinate pairs}
\label{table1}
\begin{tabular}{cc}
\hline
$\gamma$ &  $\alpha_\gamma^\ast = \beta_\gamma^\ast$ \\
\hline
$0$    & $0.158896$ \\
$0.01$ & $0.160179$ \\
$0.05$ & $0.165872$ \\
$0.10$ & $0.174668$ \\
$0.15$ & $0.186511$ \\
$0.20$ & $0.205147$ \\
$1-\pi/4$ & $1-\pi/4$ \\
\hline
\end{tabular}
\end{table}

\section{Hazard Function}

The BN hazard rate function
takes the form
\begin{align*}
h(x)
&
=
\frac{\left[\Phi\left(\frac{x-\mu}{\sigma}\right)
\right]^{\alpha-1}\,
\left[1-\Phi\left(\frac{x-\mu}{\sigma}\right)
\right]^{\beta-1}\,\phi\left(\frac{x-\mu}{\sigma}\right)}
{\sigma\,\mathrm{B}(\alpha,\beta)\,
\left[1-I_{\Phi(\frac{x-\mu}{\sigma})}(\alpha,\beta)
\right]},
\end{align*}
where
$I_x(\alpha,\beta)
=
\mathrm{B}(\alpha,\beta)^{-1}\int_0^x w^{\alpha-1}\,(1-w)^{\beta-1}
\mathrm{d}w$
denotes the incomplete beta function ratio.
At $x=\mu$, it gives
\begin{align*}
h(\mu)=
\frac{2^{2-\alpha-\beta}}{
\sqrt{2\pi}\,\sigma\,\mathrm{B}(\alpha,\beta)
\left[1-I_{\frac{1}{2}}(\alpha,\beta)\right]}.
\end{align*}
We study the asymptotic behavior of the hazard function as $x\to\infty$.
Let $h_{\text{N}}(x)$ be the hazard rate function of the normal distribution. We show that $h_{\text{N}}(x)\sim x/\sigma^2$ as $x\to\infty$. We have:
\begin{align*}
\lim_{x\to\infty}
\frac{h_{\text{N}}(x)}{x}
&=
\lim_{x\to\infty}
\frac{\phi(\frac{x-\mu}{\sigma})\frac{1}{x}}{\sigma[1-\Phi(\frac{x-\mu}{\sigma})]}.
\end{align*}
Direct evaluation of this limit gives an indeterminate form. However, using L'H\^opital rule, straightforward manipulations show that the above limit is simply $1/\sigma^2$.
We are now able to show that
$h(x)\sim (\beta/\sigma^2) x$
as $x\to\infty$.
Indeed, we have:
\begin{align*}
\lim_{x\to\infty}
\frac{h(x)}{x}
&=
\frac{1}{\mathrm{B}(\alpha,\beta)\sigma}
\lim_{x\to\infty}
\left[
\Phi\left(\frac{x-\mu}{\sigma}\right)
\right]^{\alpha-1}
\lim_{x\to\infty}
\frac{
\left[
1-\Phi\left(\frac{x-\mu}{\sigma}\right)
\right]^{\beta-1}
\phi\left(\frac{x-\mu}{\sigma}\right)\frac{1}{x}}
{1-I_{\Phi(\frac{x-\mu}{\sigma})}(\alpha,\beta)}
\\
&=
\frac{1}{\mathrm{B}(\alpha,\beta)\sigma}
\lim_{x\to\infty}
\frac{\left[
1-\Phi\left(\frac{x-\mu}{\sigma}\right)
\right]^{\beta-1}
\phi\left(\frac{x-\mu}{\sigma}\right)\frac{1}{x}}{1-I_{\Phi(\frac{x-\mu}{\sigma})}(\alpha,\beta)}.
\end{align*}
Once again another indeterminate form arises. An application of L'H\^opital rule gives
\begin{align*}
\lim_{x\to\infty}
&
\frac{
\left[
1-\Phi\left(\frac{x-\mu}{\sigma}\right)
\right]^{\beta-1}
\phi\left(\frac{x-\mu}{\sigma}\right)\frac{1}{x}}
{1-I_{\Phi(\frac{x-\mu}{\sigma})}(\alpha,\beta)}
\\
&
=
\sigma\mathrm{B}(\alpha,\beta)
\left\{
\lim_{x\to\infty}
(\beta-1)
\frac{h_{\text{N}}(x)}{x}
+
\lim_{x\to\infty}
\frac{1}{\sigma^2}
\frac{x-\mu}{x}
\right\}
\\
&
=
\mathrm{B}(\alpha,\beta)\frac{\beta}{\sigma}.
\end{align*}
Thus,
$h(x)\sim (\beta/\sigma^2)x$ as $x\to\infty$.
Combining both previous results, one can easily obtain $h(x)\sim \beta h_{\text{N}}(x)$ as $x\to\infty$.
The limit of $h(x)$ as $x\to-\infty$ is zero. However, we show that $h(x)\sim
\frac{1}{\sigma\mathrm{B}(\alpha,\beta)}\left(-\frac{x-\mu}{\sigma}\right)^{1-\alpha}
\left[\phi(\frac{x-\mu}{\sigma})\right]^\alpha$ when $x\to-\infty$. In fact, considering that
$\lim_{x\to-\infty} \Phi(\frac{x-\mu}{\sigma})=\lim_{x\to-\infty} I_{\Phi(\frac{x-\mu}{\sigma})}(\alpha,\beta)=0$, we obtain:
\begin{align*}
\lim_{x\to-\infty}
\frac{h(x)}
{
\left(-\frac{x-\mu}{\sigma}\right)^{1-\alpha}
\left[\phi(\frac{x-\mu}{\sigma})\right]^{\alpha}
}
=
\frac{1}{\sigma\mathrm{B}(\alpha,\beta)}
\lim_{x\to-\infty}
\left[
\Phi\left(\frac{x-\mu}{\sigma}\right)
\right]^{\alpha-1}
\frac{\phi\left(\frac{x-\mu}{\sigma}\right)}{
\left(
-\frac{x-\mu}{\sigma}
\right)^{1-\alpha}
\left[
\phi(\frac{x-\mu}{\sigma})
\right]^{\alpha}}.
\end{align*}
Applying L'H\^{o}pital rule, we can show that $\Phi(x)\sim (-x)^{-1}\phi(x)$. Then, we obtain
that the above limit is $1/[\sigma\mathrm{B}(\alpha,\beta)]$.

In order to relate the asymptotic behavior of
$h(x)$ and $h_{\text{N}}(x)$ as $x\to-\infty$,
let us show that
$h_{\text{N}}(x)
\sim
\phi\left(\frac{x-\mu}{\sigma}\right)
/
\sigma$
as follows
\begin{align*}
\lim_{x\to-\infty}\frac{h_{\text{N}}(x)}{\phi(\frac{x-\mu}{\sigma})}
=
\frac{1}{\sigma}
\lim_{x\to-\infty}
\frac{\phi(\frac{x-\mu}{\sigma})}
{1-\Phi(\frac{x-\mu}{\sigma})}
\frac{1}{\phi(\frac{x-\mu}{\sigma})}
=
\frac{1}{\sigma}
\lim_{x\to-\infty}
\frac{1}{1-\Phi(\frac{x-\mu}{\sigma})}=\frac{1}{\sigma}.
\end{align*}

So, the asymptotic behavior of $h(x)$ can be related to
that one of $h_{\text{N}}(x)$ when $x\to-\infty$ by:

\begin{align*}
h(x)
\sim
\frac{1}{\mathrm{B}(\alpha,\beta)\sigma^{1-\alpha}}
\left(
-\frac{x-\mu}{\sigma}
\right)^{1-\alpha}
\left[
h_{\text{N}}(x)
\right]^\alpha.
\end{align*}

\section{Moments}
We can work with the BSN distribution in generality,
since the moments of
$X=\mu+\sigma\,Z\,\sim\mathrm{BN}(\alpha,\beta,\mu,\sigma)$ follow from the
moments of $Z\sim\mathrm{BN}(\alpha,\beta,0,1)$ using
$\mathrm{E}(X^n)
=
\mathrm{E}[(\mu+\sigma Z)^n]
=
\sum_{r=0}^{n}\,\binom{n}{r}\,\mu^{n-r}\,\sigma^r\,\mathrm{E}(Z^r)$.
For $s$ and $r$ non-negative integers, let
$\tau_{s,r}=\int_{-\infty}^{\infty} x^s \phi(x)\Phi(x)^r \mathrm{d}x$
be the $(s,r)$th probability weighted moment (PWM) of the
standard normal distribution. For $\alpha$ integer and $\alpha$ real non-integer,
the $s$th moment of $Z$ can be expressed in terms of linear combinations of
these PWMs as~\cite{cordeiro2011closed}
\begin{align}
\label{moments}
\mathrm{E}(Z^s)
=
\sum_{r=0}^\infty
w_{r}(\alpha,\beta)\,\tau_{s,r+a-1}
\quad
\text{and}
\quad
\mathrm{E}(Z^s)=\sum_{i,j=0}^\infty
\sum_{r=0}^j w_{i,j,r}(\alpha,\beta)\,\tau_{s,r},
\end{align}
respectively,
where
\begin{align*}
w_i(\alpha,\beta)
=
\frac{(-1)^i\,\binom{\beta-1}{i}}{\mathrm{B}(\alpha,\beta)}
\quad
\text{and}
\quad
w_{i,j,r}(\alpha,\beta)=\frac{(-1)^{i+j+r}\,\binom{\alpha+i-1}{j}\,\binom{\beta-1}{i}\,
\binom{j}{r}}{\mathrm{B}(\alpha,\beta)}.
\end{align*}
For $s+r-l$ even, they demonstrated that
\begin{align}
\label{momentsnormal}
\begin{split}
\tau_{s,r}
=&
2^{s/2}\pi^{-(r+1/2)}
\sum_{
\substack{
l=0
\\
\text{$(s+r-l)$ even}
}
}^r
\binom{r}{l} 2^{-l}
\pi^l
\Gamma\left(\frac{s+r-l+1}{2}\right)
\times
\\
&
F_A^{(r-l)}
\left(
\frac{s+r-l+1}{2};
\frac{1}{2},\ldots,\frac{1}{2};
\frac{3}{2},\ldots,\frac{3}{2};
-1,\ldots,-1
\right),
\end{split}
\end{align}
where $F_A^{(r-l)}(\cdot)$ is the Lauricella function of type A~\cite{aarts2010lauricella}.
Expressions for terms in $\tau_{s,r}$ vanish when $s+r-l$ is odd.
Equations (\ref{moments}) and (\ref{momentsnormal}) are given as
infinite weighted sums of Lauricella functions for which numerical routines for
computation are available, for example, \mbox{Mathematica}~\cite{mathematica7}.
They extend some previously known results which are valid only for both $\alpha$ and $\beta$
integers and can be more efficient than computing the moments by writing some codes in SAS
or R. In the next section, we provide an alternative formula for the moments of the
BSN distribution based on the quantile function.

\begin{figure}
\centering
\subfigure[]{\epsfig{file=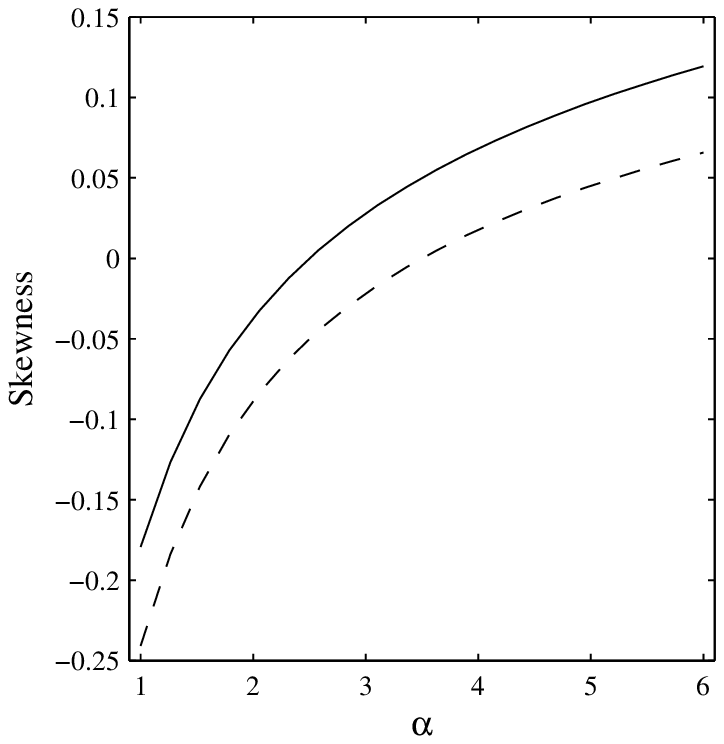}}
\subfigure[]{\epsfig{file=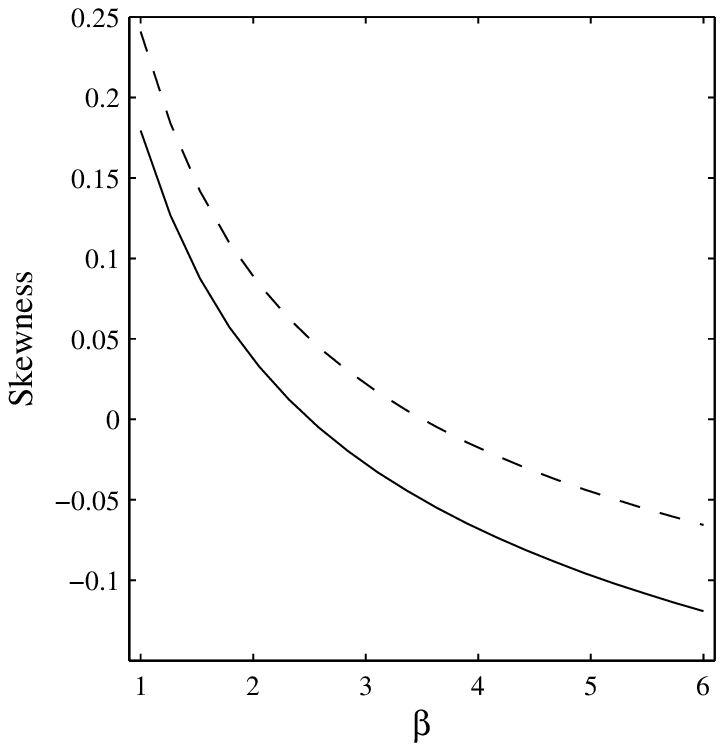}}
\caption{Plots of the BSN skewness as functions of
(a)~$\alpha$
for $\beta=2.5$ (solid curve) and $\beta=3.5$ (dashed curve),
and
(b)~$\beta$
for $\alpha=2.5$ (solid curve) and $\alpha=3.5$ (dashed curve).}
\label{fig.skew}
\end{figure}

\begin{figure}
\centering
\subfigure[]{\epsfig{file=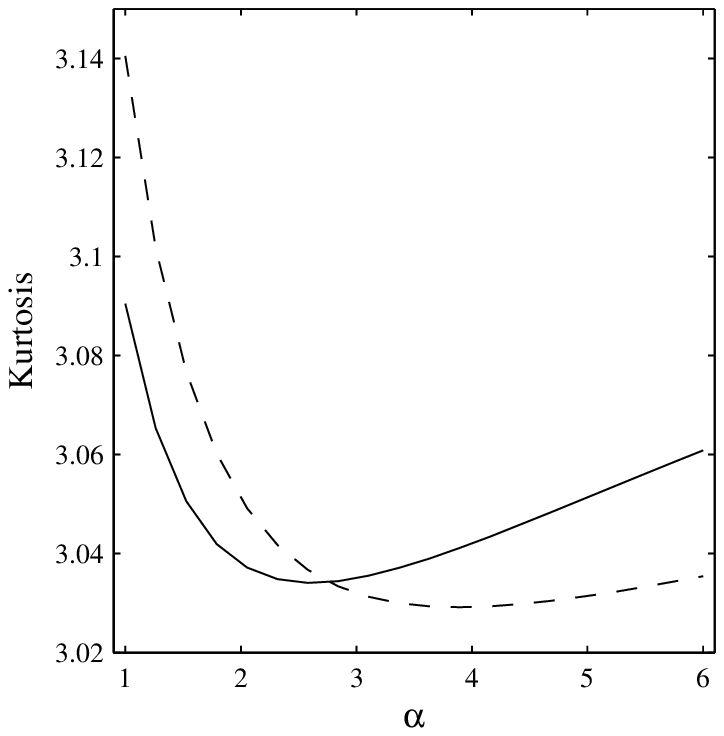}}
\subfigure[]{\epsfig{file=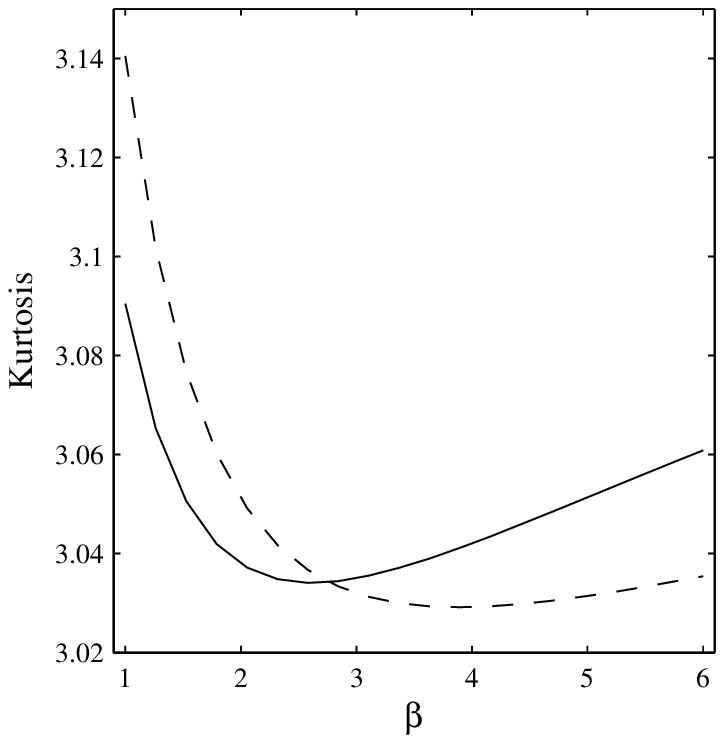}}
\caption{Plots of the BSN kurtosis as functions of
(a)~$\alpha$
for $\beta=2.5$ (solid curve) and $\beta=3.5$ (dashed curve),
and
(b)~$\beta$
for $\alpha=2.5$ (solid curve) and $\alpha=3.5$ (dashed curve).}
\label{fig.kurt}
\end{figure}

Plots of the skewness of the BSN distribution as a function of parameter $\alpha$
(for selected values of $\beta$),
and as a function of parameter $\beta$ (for selected values of $\alpha$),
are given in Figure~\ref{fig.skew}.
Figure~\ref{fig.kurt} does the same for the kurtosis of the BSN distribution.
These plots show that the BSN skewness increases
when $\alpha$ increases (for selected values of $\beta$)
and decreases when $\beta$ increases (for selected values of $\alpha$).
On the other hand, the
BSN kurtosis first decreases steadily to a minimum value and then increases
when $\alpha$ increases for fixed $\beta$ or when $\beta$ increases for fixed
$\alpha$.

\section{Quantile Function}
Quantile functions are in widespread use in general statistics and often find
representations in terms of lookup tables for key percentiles.
Without loss of
generality,
we can work with the BSN density function given by
\begin{align}\label{pdfbsn}
f(x)=\frac{1}{\mathrm{B}(\alpha,\beta)}\,\phi(x)\,\Phi(x)^{\alpha-1}\,
[1-\Phi(x)]^{\beta-1},
\end{align}
and the corresponding cumulative function reduces to
\begin{align}
\label{cdfbsn}
F(x)=I_{\Phi(x)}(\alpha,\beta)=
\frac{1}{\mathrm{B}(\alpha,\beta)}
\int_0^{\Phi(x)}
\omega^{\alpha-1}\,(1-\omega)^{\beta-1}
\mathrm{d}\omega.
\end{align}

By inverting (\ref{cdfbsn}), the BSN quantile function, say $Q(u)$, can be obtained from the
quantile functions of the standard normal and beta distributions
denoted by $Q_{\mathrm{SN}}(u)$ and $Q_{\mathrm{B}}(u)$, respectively.
We readily have $Q(u)=Q_{\mathrm{SN}}(Q_{\mathrm{B}}(u))$.
Clearly, the BSN quantile function can be calculated
from $Q_{\mathrm{SN}}(u)$ and $Q_{\mathrm{B}}(u)$ in most statistical software.

\subsection{Quantile Measures}

The effect of the shape parameters $\alpha$ and $\beta$ on the skewness and kurtosis
of the BSN distribution can be based on quantile measures. One of the earliest skewness
measures to be suggested is the Bowley skewness~\cite{kenney1962k} defined by the average
of the quartiles minus the median, divided by half the interquartile range,
namely
\begin{align*}
B=\frac{Q(3/4)+Q(1/4)-2 Q(1/2)}{Q(3/4)-Q(1/4)}.
\end{align*}
On the other hand,
the Moors kurtosis is based on octiles~\cite{moors1988kurtosis}
\begin{align*}
M=\frac{Q(7/8)-Q(5/8)+Q(3/8)-Q(1/8)}{Q(6/8)-Q(2/8)}.
\end{align*}
The measures $B$ and $M$ are less sensitive to outliers and
they exist even for distributions without moments.
Because $M$ is based on the octiles, it is not sensitive to variations
of the values in the tails or to variations of the values around
the median. Clearly, $M>0$ and there is a good agreement with
the usual kurtosis measures for some distributions.
For the normal distribution, $B=M=0$.

\begin{figure*}
\centering
\subfigure[]{\epsfig{file=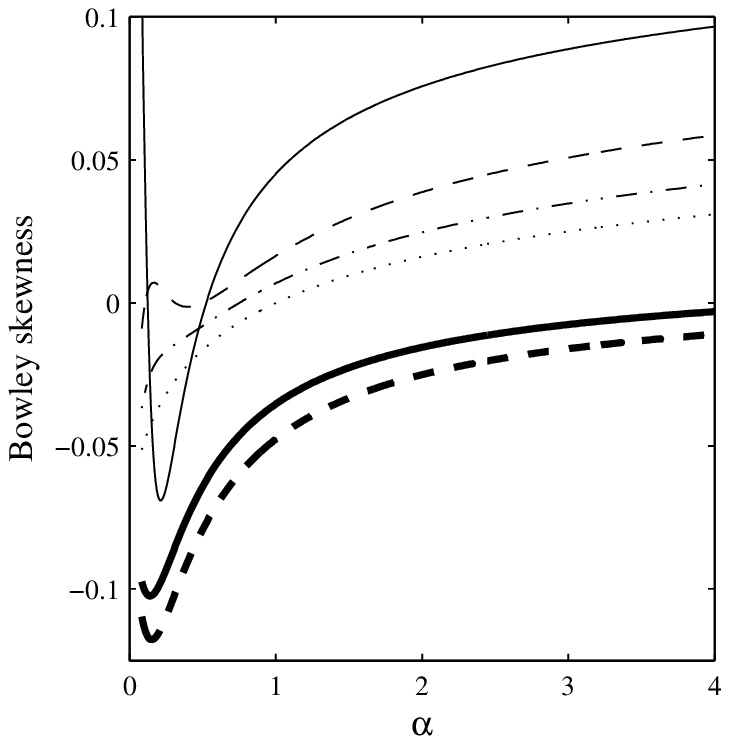}}
\subfigure[]{\epsfig{file=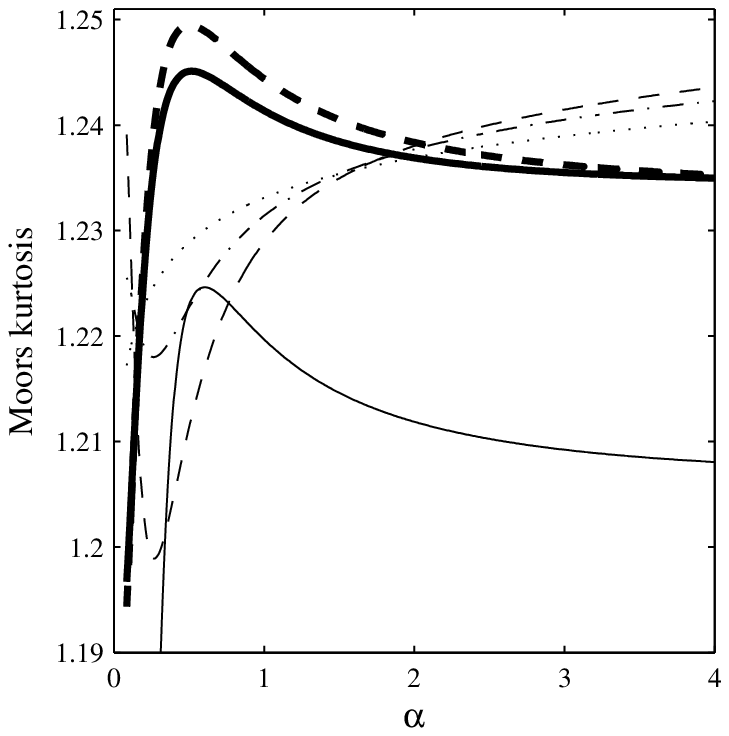}}
\caption{Bowley skewness (a) and Moors kurtosis (b) of the BSN distribution for
$0<\alpha\leq 4$
and
$\beta\in\{ 1/8, 1/2, 3/4, 1, 5, 10 \}$
(solid, dashed, dash-dotted, dotted, bold solid, and bold dashed curves,
respectively).}
\label{fig.BM-alpha}
\end{figure*}

\begin{figure*}
\centering
\subfigure[]{\epsfig{file=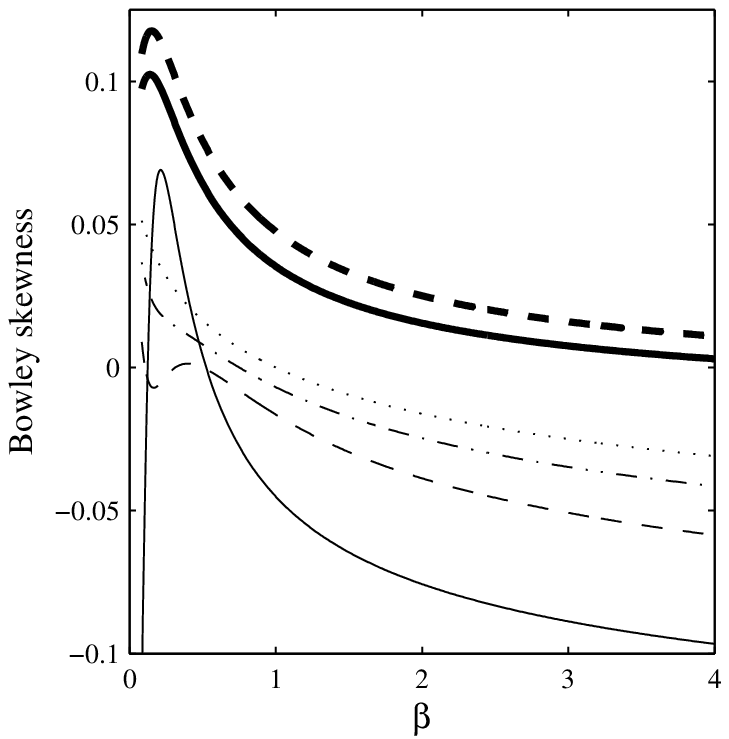}}
\subfigure[]{\epsfig{file=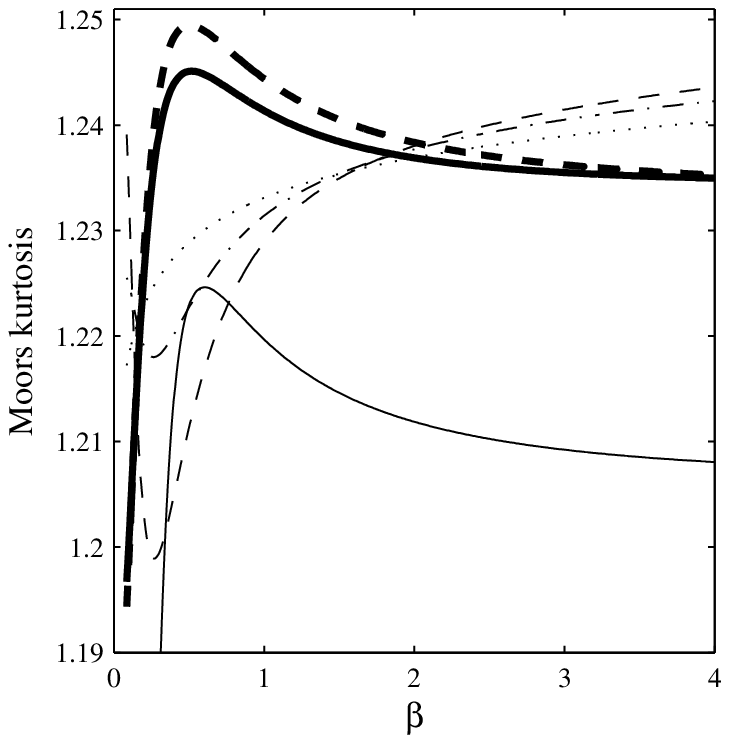}}
\caption{Bowley skewness (a) and Moors kurtosis (b) of the BSN distribution for
$0<\beta\leq 4$
and
$\alpha\in\{ 1/8, 1/2, 3/4, 1, 5, 10 \}$
(solid, dashed, dash-dotted, dotted, bold solid, and bold dashed curves,
respectively).}
\label{fig.BM-beta}
\end{figure*}

Figures~\ref{fig.BM-alpha} and \ref{fig.BM-beta}
show the measures $B$ and $M$ as functions of $\alpha$
and $\beta$ for some parameter values of the BSN distribution,
respectively.
These plots really suggest that both measures are very
sensitive on the shape parameters, thus indicating the importance of the
model~(\ref{fdpbn}).

\subsection{Power Series Expansion}\label{section-expansion-quantile}
Power series methods are at the heart of many aspects of applied
mathematics and statistics. In this section, we provide a power series
expansion for the quantile function that can be useful to obtain
mathematical properties of the BSN distribution. First, the following expansion is available for
the beta quantile function~\cite{InverseBetaRegularized2011}
\begin{align}
\label{betaquantile}
x=Q_{\mathrm{B}}(u)=I_u^{-1}(\alpha,\beta)=\sum_{i=1}^{\infty}d_i\,u^{i/\alpha}.
\end{align}
Here, $d_i=[\alpha\,\mathrm{B}(\alpha,\beta)]^{i/\alpha}\,a_i$ for $i\ge1$,
$a_0=0$, $a_1=1$, and the quantities $a_i$ for $i\ge2$ can be derived
from a cubic recurrence equation
\begin{align*}
a_{i}
=
&
\frac{1}{[i^2+(\alpha-2)i+(1-\alpha)]}
\bigg\{
(1-\delta_{i,2})\,
\sum_{r=2}^{i-1}\,a_{r}\,a_{i+1-r}\,[r(1-\alpha)(i-r)
\\
&
-
r(r-1)]
+
\sum_{r=1}^{i-1}
\sum_{s=1}^{i-r}\,
a_{r}\,a_{s}\,a_{i+1-r-s}\,
[
r(r-\alpha)+s(\alpha+\beta-2)
\\
&
\times
(i+1-r-s)
]
\bigg\}
,
\end{align*}
where $\delta_{i,2}=1$ if $i=2$ and $\delta_{i,2}=0$ if
$i\ne 2$.
In the last equation, the quadratic term only contributes for $i\ge 3$.
We have $a_2=(\beta-1)/(\alpha+1)$,
$a_3=[(\beta-1)\,(\alpha^2 +3\beta\alpha-
\alpha+5\beta-4)]/[2(\alpha+1)^2(\alpha+2)]$, and son on.
Following Steinbrecher~\cite{stein02}, the standard normal quantile function
can be expanded as $Q_{\mathrm{SN}}(u)=\sum_{k=0}^{\infty}b_k\,w^{2k+1}$,
where $w=\sqrt{2\pi}\,(u-1/2)$ and the quantities $b_k$ can be calculated
recursively from
\begin{align*}
b_{k+1}
=\frac{1}{2(2k+3)}\,\sum_{r=0}^{k}
\frac{(2r+1)\,(2k-2r+1)\,b_r\,b_{k-r}}{(r+1)\,(2r+1)}.
\end{align*}
Here, $b_0=1$, $b_1=1/6$, $b_2=7/120$, $b_3=127/7560,\ldots$
The function $Q_{\mathrm{SN}}(u)$ can be written as a power series given by
\begin{align}\label{expquantile}
Q_{\mathrm{SN}}(u)=\sum_{i=0}^{\infty}c_k\,\,(u-1/2)^{k},
\end{align}
whose $c_k$ is defined by $c_k=0$ for $k=0,2,4,\ldots$ and
$c_k=(2 \pi)^{k/2}\,b_{(k-1)/2}$ for $k=1,3,5,\ldots$
Combining (\ref{betaquantile}) and (\ref{expquantile}),
$Q(u)=Q_{\mathrm{SN}}(Q_{\mathrm{B}}(u))$ can be reduced to
\begin{eqnarray}\label{quantile1}
Q(u)=\sum_{k=0}^{\infty}c_k\,\left(\sum_{i=0}^{\infty}d_i\,\,u^{i/\alpha}\right)^k,
\end{eqnarray}
where $d_0=-1/2$ and (as before)
$d_i=\alpha^{i/\alpha}\,\mathrm{B}(\alpha,\beta)^{i/\alpha}\,a_i$
for $i\ge 1$.

By application of an equation of~\cite{gradshteyn2000table}
for a power series raised to a positive integer $k$, we have
\begin{eqnarray}\label{expnint}
\left(\sum_{i=0}^\infty d_i\,\,u^{i/\alpha}\right)^k=
\sum_{i=0}^\infty e_{k,i}\,\,u^{i/\alpha},
\end{eqnarray}
whose coefficients $e_{k,i}$ (for $i=1,2,\ldots$) can be determined
numerically from the quantities $d_i$ using the recurrence equation
(with $e_{k,0}=d_0^k$)
\begin{align}\label{coeffic}
e_{k,i}=(i\,d_0)^{-1}\sum_{m=1}^{i}[m\,(k+1)-i]\,d_m\,e_{k,i-m}.
\end{align}
The coefficient $e_{k,i}$ can be calculated from
$e_{k,0},\ldots,e_{k,i-1}$ and hence from the quantities
$d_0,\ldots,d_{i}$. Further,
it can also be given explicitly
in terms of the coefficients $d_i$,
although it is not necessary for programming numerically these expansions
in any algebraic or numerical software. From (\ref{quantile1})--(\ref{coeffic}),
we can write
\begin{align}
\label{quantilefinal}
Q(u)=\sum_{i=0}^{\infty} f_i\,\,u^{i/\alpha},
\end{align}
where $f_i=\sum_{k=0}^{\infty} c_k\,e_{k,i}$, for $i=0,1,\ldots$
The power series expansion (\ref{quantilefinal}) is the main result
of this section. Some mathematical properties of the BSN distribution
(such as ordinary moments, generating function and mean deviations)
can be directly obtained from (\ref{quantilefinal}).
For example, if $Z$ has the BSN distribution,
we can use (\ref{expnint}) and (\ref{coeffic}),
to obtain after integration
the $s$th moment of $Z$ as
\begin{align}\label{moment2}
\mathrm{E}(Z^s)=\int_{0}^{1}\,
\left(\sum_{i=0}^{\infty} f_i\,\,u^{i/\alpha}\right)^s
\mathrm{d}u=\sum_{i=0}^{\infty}\,\frac{g_{s,i}}{i/\alpha+1},
\end{align}
where $g_{s,0}=f_0^s$ and $g_{s,i}=(i\,f_0)^{-1}\,\sum_{m=1}^{i}\,[m\,(s+1)-i]\,f_m\,g_{s,i-m}$
for $s\ge 1$.

\section{Generating Function}

Let $X$ be a random variable with BSN density function and $M(t)=E[\exp(t\,X)]$
be the moment generating function (mgf) of $X$. Here, we give two
representations for $M(t)$. The first representation follows from (\ref{expnint})--(\ref{quantilefinal}) as
\begin{eqnarray}\label{mgf1}
M(t)=\sum_{s=0}^{\infty}\int_{0}^{1}\frac{t^s\,\left(\sum_{i=0}^{\infty} f_i\,u^{i/\alpha}\right)^s}{s!}\mathrm{d}u=\sum_{s,i=0}^{\infty}\frac{g_{s,i}\,\,t^s}{(i/\alpha+1)\,s!},
\end{eqnarray}
where $g_{s,i}$ was defined in Section~\ref{section-expansion-quantile}.

The second representation follows,
by expanding the binomial in (\ref{pdfbsn}),
as
\begin{eqnarray}\label{exppdf1}
f(x)=\phi(x)\,\sum_{k=0}^\infty \nu_k\,\,\Phi(x)^{k+\alpha-1},
\end{eqnarray}
where
$\nu_k=\nu_k(\alpha,\beta)=(-1)^k\,\binom{\beta-1}{k}/[\mathrm{B}(\alpha,\beta)]$.
For $\delta>0$ real non-integer, we can write
\begin{equation}\label{powerG}
\Phi(x)^{\delta}=\sum_{r=0}^\infty\,s_r(\delta)\,\Phi(x)^r,
\end{equation}
where $s_r(\delta)=\sum_{j=r}^\infty (-1)^{r+j}\,
\binom{\delta}{j}\,\binom{j}{r}$. Combining (\ref{exppdf1}) and (\ref{powerG}),
we have
\begin{eqnarray}
\label{exppdf2}
f(x)=\phi(x)\,\sum_{r=0}^\infty \pi_r\,\,\,\Phi(x)^r.
\end{eqnarray}
where $\pi_r=\sum_{k=0}^\infty \nu_k\,\,s_r(k+\alpha-1)$. From (\ref{exppdf2}), we obtain
\begin{align*}
M(-t)=\frac{1}{\sqrt 2\pi}
\sum_{r=0}^\infty \pi_r\,\,
\int_{-\infty}^{\infty}
\Phi(x)^r\,
\exp\left(-t\,x-\frac{x^2}{2}\right)
\mathrm{d}x.
\end{align*}
The standard normal cdf $\Phi(x)$ can be written as a power series expansion
$\Phi(x)=\sum_{j=0}^{\infty}a_j\,x^{j}$,
where $a_0=(1+\sqrt{2/\pi})^{-1}/2$, $a_{2j+1}=\frac{(-1)^{j}2^{-j}}{\sqrt{2\pi}(2j+1)j!}$
for $j=0,1,2\ldots$ and $a_{2j}=0$ for $j=1,2,\dots$
We can write $\Phi(x)^r=\sum_{j=0}^{\infty}c_{r,j}\,x^{j}$,
whose coefficients $c_{r,j}$ can be determined from (\ref{expnint})
and (\ref{coeffic}) by $c_{r,0}=a_0^r$ and
$c_{r,j}=(j\,a_0)^{-1}\sum_{m=1}^{i}[m\,(r+1)-j]\,a_m\,c_{r,j-m}$.
Hence,
\begin{align*}
M(-t)=\frac{1}{\sqrt 2\pi}\sum_{r,j=0}^\infty \,\pi_r\,\,c_{r,j}\,
\int_{-\infty}^{\infty}\,x^{j}\,
\exp\left(-t x-\frac{x^2}{2}\right)
\mathrm{d}x.
\end{align*}
By equation (2.3.15.8) in~\cite{prudnikov1986integrals},
the integral becomes
\begin{align*}
J(t,j)=\int_{-\infty}^{\infty}\,x^{j}\,
\exp\left(-t x-\frac{x^2}{2}\right)
\mathrm{d}x
=(-1)^j\,\sqrt{2\pi}\,\,\,\frac{\partial^j}{\partial t^j}
\left\{\exp\left(\frac{t^2}{2}\right)\right\}
\end{align*}
and thus
\begin{align}
\label{mgf2}
M(-t)=\frac{1}{\sqrt 2\pi}
\sum_{r,j=0}^\infty \,\pi_r\,\,c_{r,j}\,\,J(t,j).
\end{align}
Equations (\ref{mgf1}) and (\ref{mgf2}) are the main results of this section.

\section{Mean Deviations}

If $X$ has the BN distribution,
we can derive the mean deviations about the mean $\nu=E(X)$
and about the median $m$ from
$\delta_1=\int_{-\infty}^{\infty}\mid\!x - \nu\!\!\mid f(x) \mathrm{d}x$
and
$\delta_2=\int_{-\infty}^{\infty}\mid\! x - m\!\!\mid f(x)\mathrm{d}x$,
respectively.
The mean $\nu$ can be obtained from (\ref{moment2}) with $s=1$ and the
median $m$ is the solution of the non-linear equation
$I_{\Phi\left(\frac{m-\mu}{\sigma}\right)}(\alpha,\beta)=1/2$. The deviations
$\delta_1$ and $\delta_2$ can be expressed as
\begin{align}
\label{deltas}
\delta_1= 2\big[\nu F(\nu)-J(q)\big]
\quad\,\,\,
\text{and}
\quad\,\,\,
\delta_2=\nu-2\,J(m),
\end{align}
where $J(q)=\int_{-\infty}^{q} x\,f(x) \mathrm{d}x$.
We provide an explicit expression
for $J(q)$ based on the quantile function expansion.
From (\ref{quantilefinal}),
we have
\begin{align*}
J(q)=
\int_{0}^{I_{\Phi\left(\frac{q-\mu}{\sigma}\right)}(\alpha,\beta)}\,
Q(u)\mathrm{d}u=\sum_{i=0}^{\infty}
\frac{f_i\,\,\,I_{\Phi\left(\frac{q-\mu}{\sigma}\right)}(\alpha,\beta)^{(i/\alpha+1)}}
{i/\alpha+1}.
\end{align*}
Hence, the mean deviations follow from equations (\ref{deltas}).

\section{Shannon Entropy}

Shannon entropy is a measure of uncertainty associated with a random variable.
Consider a random variable $X\sim \mathrm{BN}(\alpha,\beta,\mu,\sigma)$. Thus,
the Shannon entropy of $X$ is given by
\begin{align*}
\mathrm{H}(X)
=&
-
\mathrm{E}\{\ln f(X)\}
=
-
\int_{-\infty}^{\infty} f(x)\ln f(x)\mathrm{d}x
\\
=
&
-\int_{\infty}^{\infty}
f(x)
\ln
\left\{
\frac{1}{\sigma\mathrm{B}(\alpha,\beta)}
\left[
\Phi\left(\frac{x-\mu}{\sigma} \right)
\right]^{\alpha-1}
\left[
1 - \Phi\left( \frac{x-\mu}{\sigma} \right)
\right]^{\beta-1}
\phi\left( \frac{x-\mu}{\sigma} \right)
\right\}
\mathrm{d}x
\\
=
&
-
\int_{-\infty}^{\infty}
f(x)
\ln
\left\{
\frac{1}{\sigma\mathrm{B}(\alpha,\beta)}
\right\}
\mathrm{d}x
-
\int_{-\infty}^{\infty}
f(x)
(\alpha-1)
\ln
\left\{
\Phi\left( \frac{x-\mu}{\sigma}\right)
\right\}
\mathrm{d}x
\\
&
-
\int_{-\infty}^{\infty}f(x) (\beta-1)
\ln
\left\{
1-\Phi\left( \frac{x-\mu}{\sigma}\right)
\right\}
\mathrm{d}x
-
\int_{-\infty}^{\infty}f(x)
\ln
\left\{
\phi\left( \frac{x-\mu}{\sigma} \right)
\right\}
\mathrm{d}x.
\end{align*}
The first integral is simply equal to $-\ln[\sigma\mathrm{B}(\alpha,\beta)]$.
The second and third integrals can be expressed in terms of the beta function.
Indeed, using Taylor series expansions for the logarithmic function,
we obtain
\begin{align*}
&\ln\left\{
\Phi\left( \frac{x-\mu}{\sigma}\right)\right\}=
\sum_{n=1}^{\infty}\frac{(-1)^n}{n}\,\left[
1-\Phi\left( \frac{x-\mu}{\sigma}\right)\right]^n
\end{align*}
and
\begin{align*}
\ln\left\{1-\Phi\left( \frac{x-\mu}{\sigma}\right)\right\}
=\sum_{n=1}^{\infty}\frac{(-1)^n}{n}\,
\Phi\left( \frac{x-\mu}{\sigma}\right)^n.
\end{align*}
As a consequence,
the previous integrals reduce to
$$
\int_{-\infty}^{\infty}
f(x)(\alpha-1)
\ln\left\{\Phi\left(\frac{x-\mu}{\sigma}\right)
\right\}\mathrm{d}x=
\frac{1-\alpha}{\mathrm{B}(\alpha,\beta)}
\sum_{n=1}^{\infty}\frac{\mathrm{B}(\alpha,n+\beta)}{n}
$$
and
$$
\int_{-\infty}^{\infty}
f(x)(\beta-1)
\ln\left\{1-\Phi\left(\frac{x-\mu}{\sigma}\right)
\right\}
\mathrm{d}x
=\frac{1-\beta}{\mathrm{B}(\alpha,\beta)}
\sum_{n=1}^{\infty}\frac{\mathrm{B}(n+\alpha,\beta)}{n}.
$$
Now, we can write the last integral as
\begin{align*}
\int_{-\infty}^{\infty}f(x)\,\ln\left\{
\phi\left( \frac{x-\mu}{\sigma} \right)
\right\}
\mathrm{d}x
&
=
\int_{-\infty}^{\infty}
f(x)
\left[
\ln
\left(
\frac{1}{\sqrt{2\pi}}
\right)
+
\ln\left\{
\exp\left(-\frac{(x-\mu)^2}{2\sigma^2}\right)
\right\}
\right]
\mathrm{d}x
\\
&
=
-\ln\left(\sqrt{2\pi}\right)
-
\frac{1}{2\sigma^2}
\left[
\mathrm{E}(X^2)-2\mu \mathrm{E}(X) + \mu^2
\right]
,
\end{align*}
where $\mathrm{E}(X^2)$ and $\mathrm{E}(X)$ follow from
(\ref{moment2}). These moments are also
given in~\cite{gupta2004moments}. Finally, we conclude that

\begin{align*}
\mathrm{H}(X)
=&
\ln\left\{
\sqrt{2\pi} \sigma \mathrm{B}(\alpha,\beta)
\right\}
+
\frac{1}{2\sigma^2}
\left[
\mathrm{E}(X^2)-2\mu \mathrm{E}(X) + \mu^2
\right]
\\
&
+
\frac{1}{\mathrm{B}(\alpha,\beta)}
\sum_{n=1}^{\infty}
\frac{1}{n}
\left[
(\alpha-1)\mathrm{B}(\alpha,n+\beta)+(\beta-1)\mathrm{B}(n+\alpha,\beta)
\right].
\end{align*}

{\small
\singlespacing
\bibliographystyle{siam}
\bibliography{bnc}
}

\end{document}

%% file: figure4a.pstex_t
\begin{picture}(0,0)%
\includegraphics{figure4a.pstex}%
\end{picture}%
\setlength{\unitlength}{4144sp}%
\begingroup\makeatletter\ifx\SetFigFont\undefined%
\gdef\SetFigFont#1#2#3#4#5{%
  \reset@font\fontsize{#1}{#2pt}%
  \fontfamily{#3}\fontseries{#4}\fontshape{#5}%
  \selectfont}%
\fi\endgroup%
\begin{picture}(3397,3461)(451,-2802)
\put(2566,-1366){\makebox(0,0)[lb]{\smash{{\SetFigFont{9}{10.8}{\rmdefault}{\mddefault}{\updefault}{\color[rgb]{0,0,0}Increasing $\beta$}%
}}}}
\end{picture}%

%% file: figure4b.pstex_t
\begin{picture}(0,0)%
\includegraphics{figure4b.pstex}%
\end{picture}%
\setlength{\unitlength}{4144sp}%
\begingroup\makeatletter\ifx\SetFigFont\undefined%
\gdef\SetFigFont#1#2#3#4#5{%
  \reset@font\fontsize{#1}{#2pt}%
  \fontfamily{#3}\fontseries{#4}\fontshape{#5}%
  \selectfont}%
\fi\endgroup%
\begin{picture}(3397,3461)(451,-2802)
\put(2566,-556){\makebox(0,0)[lb]{\smash{{\SetFigFont{9}{10.8}{\rmdefault}{\mddefault}{\updefault}{\color[rgb]{0,0,0}Increasing $\alpha$}%
}}}}
\end{picture}%